\documentclass[letter,11pt]{article}

\usepackage{amsthm} 
\usepackage{amsfonts, mathrsfs,amssymb}
\usepackage{amsmath}
\usepackage[numbers, sort&compress]{natbib}
\usepackage{graphicx}
\usepackage{vmargin, enumerate}
\usepackage{paralist}
\usepackage{times}
\usepackage[pdftex]{hyperref}

\newtheorem{thm}{Theorem}

\newtheorem{lem}[thm]{Lemma}
\newtheorem{prop}[thm]{Proposition}

\setmarginsrb{1in}{1in}{1in}{.7in}{0pt}{0pt}{.1in}{.3in}

\newcommand{\N}{\ensuremath{\mathbb N}}

\newcommand{\Do} {\mathcal{D}[0,1]}

\newcommand{\Drn}{\mathcal{D}_{r_n}[0,1]}

\newcommand{\p}[1]{{\mathbf P}\left(#1\right)}

\newcommand{\event}[1]{\left[\, #1 \,\right]}

\newcommand{\pran}[1]{\left( #1 \right)}

\newcommand{\Ec}[1]{\ensuremath{\mathbf{E} [#1]}}

\newcommand{\CE}[2]{\mathbf{E}\event{\left. #1 \; \right| \; #2 \;}}

\newcommand{\eqdist}{\ensuremath{\stackrel{d}{=}}}

\newcommand{\E}[1]{\ensuremath{\mathbf{E} \left[#1 \right]}}
\newcommand{\V}[1]{\ensuremath{\mathbf{Var} \left(#1 \right)}}
\newcommand{\Vc}[1]{\ensuremath{\mathbf{Var} (#1)}}
\newcommand{\Prob}[1]{\ensuremath{\mathbf{P} \left(#1 \right)}}
\newcommand{\I}[1]{\ensuremath{\mathbf{1}_{ \{ #1 \} }}}
\newcommand{\fl}[1]{\ensuremath{\lfloor #1 \rfloor}}

\newcommand{\convdist}{\ensuremath{\stackrel{d}{\rightarrow}}}

\newcommand{\sq}[1]{\left[#1\right]}

\newcommand{\vol}{\ensuremath{\textup{Leb}}}

\newcommand{\Be}{\textup{B}}

\hypersetup{
    bookmarks=true,         
    unicode=false,          
    pdftoolbar=true,        
    pdfmenubar=true,        
    pdffitwindow=true,      
    pdftitle={My title},    
    pdfauthor={Author},     
    pdfsubject={Subject},   
    pdfnewwindow=true,      
    pdfkeywords={keywords}, 
    colorlinks=true,       
    linkcolor=blue,          
    citecolor=blue,        
    filecolor=blue,      
    urlcolor=blue           
 }

\begin{document}

\title{\bf Partial match queries in random quadtrees}
\author{Nicolas Broutin \\ nicolas.broutin@inria.fr \\ Projet Algorithms \\
INRIA Rocquencourt \\
78153 Le Chesnay\\ France  \and Ralph Neininger and  Henning Sulzbach  \\  \{neiningr, sulzbach\}@math.uni-frankfurt.de\\ Institute for Mathematics (FB 12) \\
J.W.~Goethe University\\ 60054 Frankfurt am Main\\ Germany} 

\maketitle

\begin{abstract}
We consider the problem of recovering items matching a partially specified pattern in multidimensional trees (quad trees and k-d trees). We assume the traditional model where the data consist of independent and uniform points in the unit square. For this model, in a structure on $n$ points, it is known that the number of nodes $C_n(\xi)$ to visit in order to report the items matching an independent and uniformly on $[0,1]$ random query $\xi$ satisfies $\Ec{C_n(\xi)}\sim \kappa n^{\beta}$, where $\kappa$ and $\beta$ are explicit constants. We develop an approach based on the analysis of the cost $C_n(x)$ of any fixed query $x\in [0,1]$, and give precise estimates for the variance and limit distribution of the cost $C_n(x)$. Our results permit to describe a limit process for the costs $C_n(x)$ as $x$ varies in $[0,1]$; one of the consequences is that $\Ec{\max_{x\in [0,1]} C_n(x)} \sim \gamma n^\beta$ ; this settles a question of Devroye [Pers.\ Comm., 2000].
\end{abstract}

\newcommand{\noi}{\noindent}

\section{Introduction}\label{sec:intro}


Multidimensional databases arise in a number of contexts such as computer graphics, management of geographical data or statistical analysis. The question of retrieving the data matching a specified pattern is then of course of prime importance. If the pattern specifies all the data fields, the query can generally be answered in logarithmic time, and a great deal of precise analyses are available in this case \cite{FlLa1994,FlLaLaSa1995,Knuth1998,Mahmoud1992a,FlSe2009}. We will be interested in the case when the pattern only constrains some of the data fields; we then talk of a \emph{partial match query}. 

The first investigations about partial match queries by \citet{Rivest1976} were based on digital structures. In a comparison-based setting, a few general purpose data structures generalizing binary search trees permit to answer partial match queries, namely the quadtree \cite{FiBe1974}, the $k$-d tree \cite{Bentley1975} and the relaxed $k$-d tree \cite{DuEsMa1998}. 
Aside of the interest that one might have in partial match for itself, there are numerous reasons that justify the precise quantification of the cost of such general search queries in comparison-based data structures. The high dimesional trees are indeed a data structure of choice for applications that range from collision detection in motion planning to mesh generation that takes advantage of the adaptive partition of space that is produced \cite{YeSh1983a,hole1988}. For general references on multidimensional data structures and more details about their various applications, see the series of monographs by \citet{Samet1990a,Samet1990,Samet2006}. The cost of partial match queries also appears in (hence influences) the complexity of a number of other geometrical search questions such as range search \cite{DuMa2002a} or rank selection \cite{DuJiMa2010}. 

In spite of its importance, the complexity results about partial match queries are not as precise as one could expect. In this paper, we provide novel analyses of the costs of partial match queries in some of the most important two dimensional data structures. Most of the document will focus on the special case of quadtrees ; in a final section, we discuss the case of $k$-d tree \cite{Bentley1975} and relaxed $k$-d trees \cite{DuEsMa1998}.

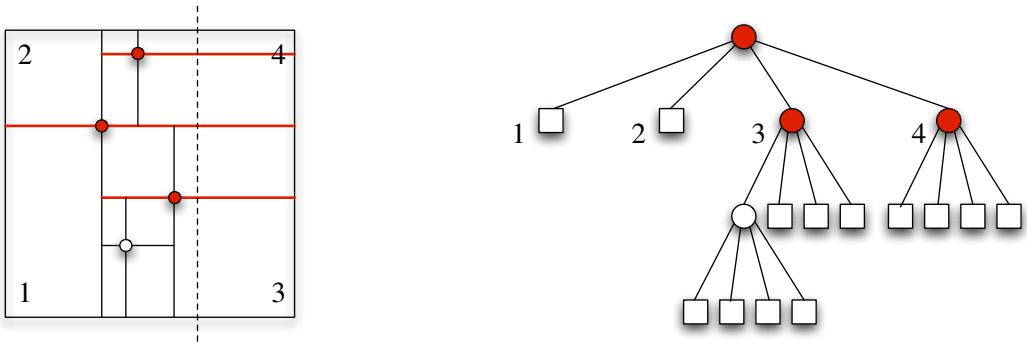
\begin{figure}
\begin{picture}(300,120)
\put(30,0){\includegraphics[scale=.5]{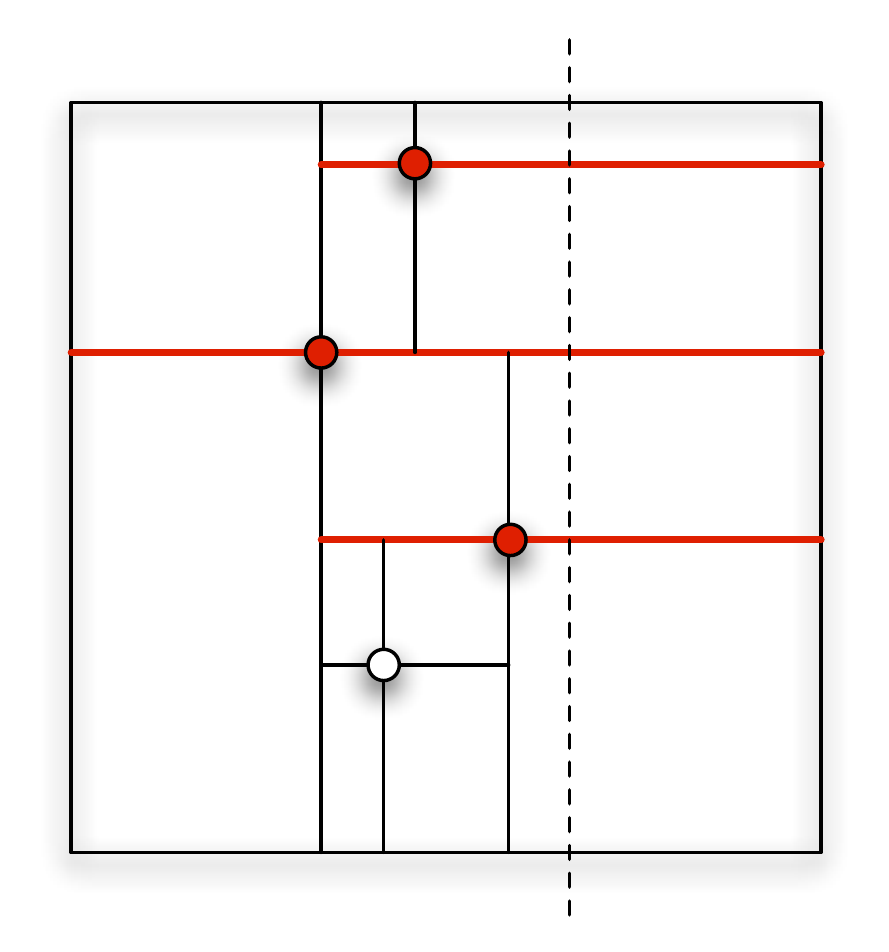}}
\put(45,20){1}
\put(45,110){2}
\put(140,20){3}
\put(140,110){4}
\put(230,0){\includegraphics[scale=.5]{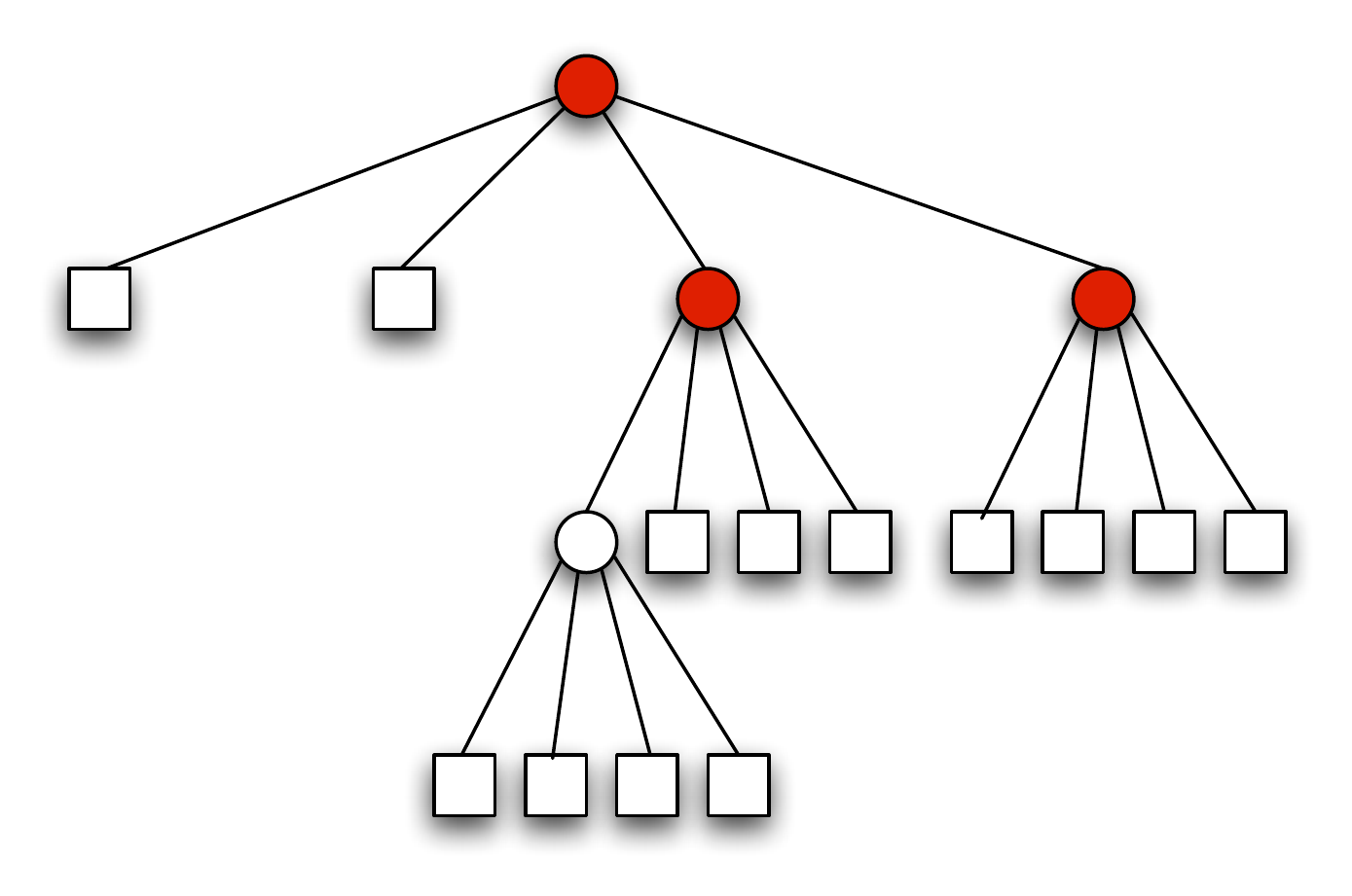}}
\put(230,80){1}
\put(275,80){2}
\put(320,80){3}
\put(380,80){4}
\end{picture}
\caption{\label{fig:quadtree}\small An example of a (point) quadtree: on the left the partition of the unit square induced by the tree data structure on the right (the children are ordered according to the numbering of the regions on the left). Answering the partial match query materialized by the dashed line on the left requires to visit the points/nodes coloured in red. Note that each one of the visited nodes correspond to a horizontal line that is crossed by the query.}
\end{figure}

\medskip
\noi\textsc{Quad trees and multidimensional search.}
The quadtree \cite{FiBe1974} allows to manage multidimensional data by extending the divide-and-conquer approach of the binary search tree. 
Consider the point sequence $p_1,p_2,\dots, p_n\in [0,1]^2$. As we build the tree, regions of the unit square are associated to the nodes where the points are stored. Initially, the root is associated with the region $[0,1]^2$ and the data structure is empty. The first point $p_1$ is stored at the root, and divides the unit square into four regions $Q_1, \dots, Q_4$. Each region is assigned to a child of the root. More generally, when $i$ points have already been inserted, we have a set of $1+3^i$ (lower-level) regions that cover the unit square. The point $p_{i+1}$ is stored in the node (say $u$) that corresponds to the region it falls in, divides it into four new regions that are assigned to the children of $u$. See Figure~\ref{fig:quadtree}. 

\medskip
\noi\textsc{Analysis of partial match retrieval.}
For the analysis, we will focus on the model of \emph{random quadtrees}, where the data points are uniformly distributed in the unit square. In the present case, the data are just points, and the problem of partial match retrieval consists in reporting all the data with one of the coordinates (say the first) being $s\in [0,1]$.
It is a simple observation that the number of nodes of the tree visited when performing the search is precisely $C_n(s)$, the number of 
regions in the quadtree that insersect a vertical line at $s$. The first analysis of partial match in quadtrees is due to \citet{FlGoPuRo1993} (after the pioneering work of \citet{FlPu1986} in the case of $k$-d trees). They studied the singularities of a differential system for the generating 
functions of partial match cost to prove that, for a random query $\xi$, being independent of the tree and uniformly distributed on $[0,1]$,
\begin{equation} \label{def:beta}
 \Ec{C_n(\xi)}\sim \kappa \: n^{\beta}\qquad \text{where}\qquad \kappa = \frac{\Gamma(2\beta+2)}{2\Gamma(\beta+1)^3}, \quad \beta=\frac{\sqrt{17}-3}2,
\end{equation}
and $\Gamma(x)$ denotes the Gamma function $\Gamma(x)=\int_0^\infty t^{x-1} e^{-t} dt$. This has since been strengthened by \citet{ChHw2003}, who provided the order of the error term (together with the values of the leading constant in all dimensions). The most precise result is (6.2) there, saying that
\begin{equation} \label{hwangquad}
\Ec{C_n(\xi)} = \kappa \: n^{\beta} - 1 + O(n^{\beta -1}).
\end{equation}

To gain a refined understanding of the cost beyond the level of expectations we pursue two directions. First, to justify that the expected value is a reasonable estimate of the cost, one would like a guarantee that the cost of partial match retrieval are actually close to their mean. However, deriving higher moments turns out to be more subtle than it seems. In particular, when the query line is random (like in the uniform case) although the four subtrees at the root really are independent given their sizes, the contributions of the two subtrees that do hit the query line are \emph{dependent}! The relative location of the query line inside these two subtrees, is again uniform, but unfortunately it is same in both regions. This issue has not yet been addressed appropriately, and there is currently no result on the variance of or higher moments for $C_n(\xi)$.

The second issue lies in the very definition of the cost measure: even if the data follow some distribution (here uniform), should one really assume that the query also satisfies this distribution? In other words, should we focus on $C_n(\xi)$? Maybe not. But then, what distribution should one use for the query line?

One possible approach to overcome both problems is to consider the query line to be fixed and to study $C_n(s)$ for $s\in[0,1]$. This raises another problem: even if $s$ is fixed at the top level, as the search is performed, the \emph{relative} location of the the queries in the recursive calls varies from a node to another! Thus, in following this approach, one is led to consider the entire process $C_n(s), s\in [0,1]$ ; this is the method we use here.  

Recently \citet{CuJo2010} obtained some results in this direction. They proved that for every fixed $s\in(0,1)$,
\begin{align} \label{const:CJ}
\Ec{C_n(s)}\sim K_1(s(1-s))^{\beta/2} n^{\beta}, \qquad K_1 
= \frac{\Gamma(2 \beta +2) \Gamma(\beta + 2)}{2 \Gamma(\beta + 1)^3 \Gamma \left (\beta/2 + 1 \right)^2}.
\end{align} 
On the other hand, \citet{FlGoPuRo1993,FlLaLaSa1995} prove that, along the edge one has $\Ec{C_n(0)}= \Theta (n^{\sqrt{2}-1})=o(n^{\beta})$ (see also \cite{CuJo2010}). The behaviour about the $x$-coordinate $U$ of the first data point certainly resembles that along the edge, so that one has $\Ec{C_n(U)}=o(n^{\beta})$. It suggests that $C_n(s)$ should not be concentrated around its mean, 
and that $n^{-\beta}C_n(s)$ should converge to a non-trivial random variable as $n\to\infty$. This random variable would of course carry much information about the asymptotic properties of the cost of partial match queries in quadtrees. Below, we identify these limit random variables and obtain refined asymptotic information on the complexity of partial match queries in quadtrees from them.


\section{Main results and implications}
\begin{figure}
	\includegraphics[scale=.6]{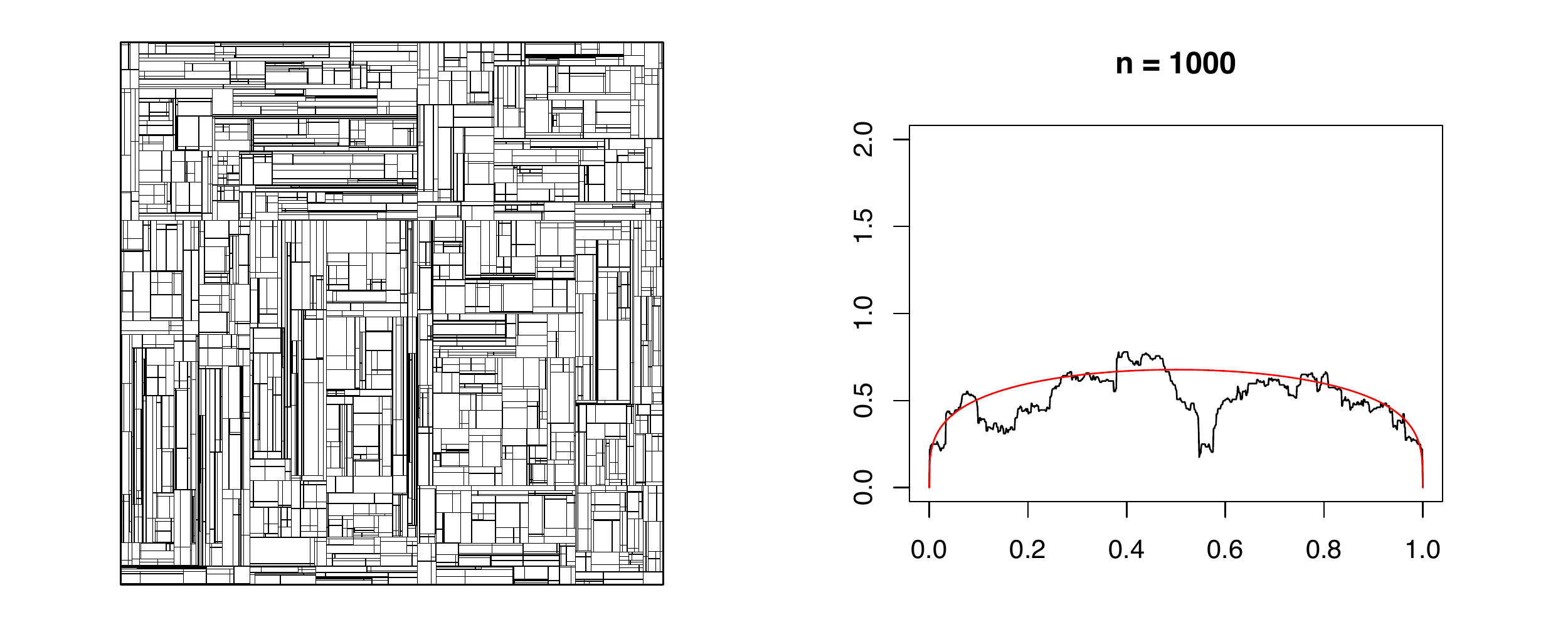}
	\caption{\label{fig:limit_process}A random quadtree on $1000$ points and the corresponding partial match process on the right ; in red we have shown the expected value.}
\end{figure}

Our main contribution is to prove  the following convergence result:
\begin{thm}\label{thm:process}Let $C_n(s)$ be the cost of a partial match query at a fixed line $s$ in a random quadtree. Then, there exists a random continuous function $Z$ such that, as $n\to\infty$, 
\begin{equation}\label{eq:process}
\pran{\frac{C_n(s)}{K_1 n^{\beta}},s\in [0,1]} \convdist (Z(s), s\in [0,1]).
\end{equation} 
This convergence in distribution holds in the Banach space $(\Do, \| \cdot \|)$ of right-continuous functions with left limits (c\`adl\`ag) equipped with the supremum norm  defined by $\|f\|=\sup_{s\in[0,1]}|f(s)|$.
\end{thm}
Note that the convergence in \eqref{eq:process} above is stronger than the convergence in distribution of the finite dimensional marginals 
$$\pran{ \frac {C_n(s_1)}{K_1 n^{\beta}}, \frac {C_n(s_2)}{K_1 n^{\beta}},\dots, \frac {C_n(s_k)}{K_1 n^{\beta}}}\convdist (Z(s_1),Z(s_2), \dots, Z(s_k))$$
as $n\to\infty$, for any natural number $k$ and points $s_1,s_2,\dots, s_k\in [0,1]$ \cite[see, e.g.,][]{Billingsley1999}.
Theorem~\ref{thm:process} has a myriad of consequences in terms of estimates of the costs of partial match queries in random quadtrees. Of course, Theorem~\ref{thm:process} would be of less practical interest if we could not characterize the distribution of the random function $Z$ (see Figure~\ref{fig:limit_process} for a simulation):
\begin{prop}
The distribution of the random function $Z$ in \eqref{eq:process} is a fixed point of the following recursive functional equation
\begin{align}\label{eq:limit_process_a}
Z(s)\eqdist 
&\I{s<U} \sq{(UV)^\beta Z^{(1)}\pran{\frac s U}+(U(1-V))^\beta Z^{(2)}\pran{\frac s U}}\nonumber\\
&+\I{s\ge U} \sq{((1-U)V)^\beta Z^{(3)}
\pran{\frac {s-U}{1-U}}+((1-U)(1-V))^\beta Z^{(4)}\pran{\frac{s-U}{1-U}}},
\end{align}
where $U$ and $V$ are independent $[0,1]$-uniform random variables and $Z^{(i)}$, $i=1,\dots, 4$ are independent copies of the process $Z$, which are also independent of $U$ and $V$. Furthermore, $Z$ in \eqref{eq:process} is the only solution of \eqref{eq:limit_process_a} such that $\Ec{Z(s)}= (s(1-s))^{\beta/2}$  for all $s\in [0,1]$ and $\Ec{\|Z\|^2}<\infty$.
\end{prop}

This is indeed relevant since the convergence that implies Theorem~\ref{thm:process} is strong enough to guarantee convergence of the variance of the costs of partial match queries. The following theorem for uniform queries $\xi$ is the direct extension of the pioneering work of \citet{FlPu1986,FlGoPuRo1993} for the cost of partial match queries at a uniform line in random multidimensional trees.
\begin{thm} \label{thm:moments}
If $\xi$ is uniformly distributed on $[0,1]$, independent of $(C_n)$ and $Z$, then 
$$\frac{C_n(\xi)}{K_1 n^{\beta}} \convdist Z(\xi),$$ 
in distribution with convergence of the first two moments. In particular 
$$\V{C_n(\xi)} \sim K_4 n^{2 \beta} \qquad \text{where} \qquad K_4 := K_1^2 \cdot \Vc{Z(\xi)} \approx 0.447363034.$$
\end{thm}
In particular, Theorem~\ref{thm:moments} identifies the asymptotic order of $\Vc{C_n(\xi)}$   which is to be compared with studies that neglected the dependence between the contributions of the subtrees mentioned above \cite{MaPaPr2001,NeRu2001,Neininger2000}. 
We also have an asymptotic for the variance of the cost at a fixed query:

\begin{thm} \label{thm:variance_fixed}
We have for all $s\in(0,1)$, as $n\to \infty$,
\begin{align} \label{const:varfix}
 \V{C_n(s)}\sim \left(2 \emph{B}(\beta +1, \beta +1) \frac{2\beta+1}{3(1-\beta)} -1\right) (s(1-s))^\beta n^{2\beta}.
\end{align}
Here, $\emph{B}(a,b):=\int_0^1 x^{a-1} (1-x)^{b-1} \,dx$ denotes the Eulerian beta integral ($a,b>0$).
\end{thm}

Some of the most striking consequence concerns the cost of the \emph{worst query} in a random plane quadtree. Note in particular that the supremum does not induce any extra logarithmic terms in the asymptotic cost.  
\begin{thm}\label{thm:supremum}Let $S_n=\sup_{s\in[0,1]} C_n(s)$. Then, as $n\to\infty$,
$$n^{-\beta} S_n \convdist S \eqdist \sup_{s\in [0,1]} Z(s)\qquad \text{and} \qquad \Ec{S_n} \sim n^{\beta} \Ec{S}, \qquad \Vc{S_n} \sim n^{2\beta} \Vc{S}.$$
\end{thm}

Finally we note  that the one-dimension marginals of the limit process $(Z(s), s\in [0,1])$ are all the same up to a multiplicative constant. 
\begin{thm}
There is a random variable $Z\ge 0$  such that for all $s\in[0,1]$,
\begin{align}\label{1dimrv}
Z(s)\stackrel{d}{=}  (s(1-s))^{\beta/2} Z.
\end{align}
The distribution of $Z$ is characterized by its moments $c_m:=\E{Z^m}$, $m\in \N$. They are given by $c_1=1$ and the recurrence 
\begin{align*}
c_m = \frac{2(\beta m +1)}{(m-1) \left(m+1 - \frac 3 2 \beta m \right)} \sum_{\ell = 1}^{m-1} {m \choose \ell} \emph{B}(\beta \ell +1, \beta (m-\ell)+1) c_\ell c_{m-\ell},  \quad m\ge 2.
\end{align*}
\end{thm}

\medskip
\noindent\textsc{Plan of the paper.} Our approach requires to work with random functions; as one might expect, proving convergence in a space of functions involves a fair amount of unavoidable technicalities. Here, we try to keep the discussion at a rather high level, to avoid diluting the main ideas in an ocean of intricate details. In Section~\ref{sec:contraction}, we give an overview of our main tool, the contraction method. 
In Section~\ref{sec:moments_sup}, we identify the variance and the supremum of the limit process $Z$, and deduce the large $n$ asymptotics for $C_n(s)$ in Theorems~\ref{thm:moments} and~\ref{thm:supremum}.

\section{Contraction method: from the real line to functional spaces} \label{sec:contraction}
\subsection{Overview}
The aim of this section is give an overview of the method we employ to prove Theorem~\ref{thm:process}. The idea is very natural and relies on a contraction argument in a certain space of probability distributions. In the context of the analysis of performance of algorithms, the method was first employed by \citet{Roesler1991a} who proved convergence in distribution for the rescaled total cost of the randomized version of quicksort. The method was then further developed by \citet{RaRu1995}, \citet{Roesler1992}, and later on in \cite{Ro01, Ne01, NeRu04, NeRu04b, DrJaNe, EiRu07} and has permitted numerous analyses in distribution for random discrete structures. 

So far, the method has mostly been used to analyze random variables taking real values, though a few applications on functions spaces have been made, see \cite{DrJaNe, EiRu07, GR09}. 
 Here we are interested in the function space ${\cal D}[0,1]$ with the uniform topology, but the main idea persists: (1) devise a recursive equation for the quantity of interest (here the process$(C_n(s), s\in[0,1])$), and (2) prove that a properly rescaled version of the quantity converges to a fixed point of a certain map related to the recursive equation ; (3) if the map is a contraction in a certain metric space, then a fixed point is unique and may be obtained by iteration. We now move on to the first step of this program.

Write $I_1^{(n)},\dots, I_4^{(n)}$ for the number of points falling in the four regions created by the point stored at the root.
Then, given the coordinates of the first data point $(U,V)$, we have, cf.~Figure~\ref{fig:quadtree},
$$(I_1^{(n)},\dots, I_4^{(n)}) \eqdist \text{Mult}(n-1; UV,U(1-V), (1-U)(1-V), (1-U)V).$$
Observe that, for the cost inside a subregion, what matters is the location of the query line \emph{relative} to the region. Thus a decomposition at the root yields the following recursive relation, for any $n\ge 1$,
\begin{align}\label{eq:Cn_rec}
C_n(s) \eqdist 1 &+ \I{s<U} \sq{C^{(1)}_{I_1^{(n)}}\!\!\pran{\frac s U}\!+ C^{(2)}_{I_2^{(n)}}\!\!\pran{\frac s U}}\!\!+\I{s\ge U} \sq{C^{(3)}_{I_3^{(n)}}\!\!\pran{\frac{1-s}{1-U}}\!
 + C^{(4)}_{I_4^{(n)}}\!\!\pran{\frac {1-s}{1-U}}},
\end{align}
%
%
where $U,I_1^{(n)},\dots, I_4^{(n)}$ are the quantities already introduced and $(C^{(1)}_k), \dots, (C^{(4)}_k)$ are independent copies of the sequence $(C_k, k\ge 0)$, independent of $(U,V,I_1^{(n)},\dots, I_4^{(n)})$. We stress that this equation does not only hold true pointwise for fixed $s$ but also as c\'adl\'ag functions on the unit interval. The relation in \eqref{eq:Cn_rec} is the fundamental equation for us.

Letting $n\to\infty$ (formally) in \eqref{eq:Cn_rec} suggests that, if $n^{-\beta} C_n(s)$ does converge to a random variable $Z(s)$ in a sense to be precised, then the distribution of the process $(Z(s), 0\le s\le 1)$ should satisfy the following fixed point equation
\begin{align}\label{eq:limit_process}
Z(s)\eqdist 
&\I{s<U} \sq{(UV)^\beta Z^{(1)}\pran{\frac s U}+(U(1-V))^\beta Z^{(2)}\pran{\frac s U}}\nonumber\\
&+\I{s\ge U} \sq{((1-U)V)^\beta Z^{(3)}
\pran{\frac {s-U}{1-U}}+((1-U)(1-V))^\beta Z^{(4)}\pran{\frac{s-U}{1-U}}},
\end{align}
where $U$ and $V$ are independent $[0,1]$-uniform random variables and $Z^{(i)}$, $i=1,\dots, 4$ are independent copies of the process $Z$, which are also independent of $U$ and $V$.

The last step leading to the fixed point equation \eqref{eq:limit_process} needs now to be made rigorous. This is at this point that the contraction method enters the game. The  distribution of a solution to our fixed-point equation (\ref{eq:limit_process}) lies  in the  set of probability measures on the Banach space $(\mathcal D[0,1], \|\cdot\|)$, which is the set  we have to endow with a metric.
The recursive equation (\ref{eq:Cn_rec}) is an example for the following, more general setting of random additive recurrences: Let $(X_n)$ be $\Do$-valued random variables with
\begin{equation} \label{eq:rec_gen} 
X_n \stackrel{d}{=} \sum_{r=1}^K A_r^{(n)}\!\!\left( X_{I_r^{(n)}}^{(r)} \right) + b^{(n)}, \quad n \geq 1,
\end{equation}
where $(A_1^{(n)}, \ldots, A_K^{(n)})$ are random linear and continuous operators on $\Do$, $b^{(n)}$ is a $\Do$-valued random variable, $I_1^{(n)}, \ldots, I_K^{(n)}$ are random integers between $0$ and $n-1$ and $(X_n^{(1)}), \ldots, (X_n^{(K)})$ are distributed like $(X_n)$. Moreover $(A_1^{(n)}, \ldots, A_K^{(n)}, b^{(n)}, I_1^{(n)}, \ldots, I_K^{(n)}), (X_n^{(1)}), \ldots, (X_n^{(K)})$ are independent. 

To establish Theorem \ref{thm:process} as a special case of this setting we use Proposition \ref{prop:cont} below.  Proposition~\ref{prop:cont}  is covered by the forthcoming paper \cite{NeSu}. We first state conditions needed to deal with the general recurrence (\ref{eq:rec_gen}); we will then justify that it can indeed be used in the case of cost of partial match queries. Consider the following assumptions, where, for a random linear operator $A$ we write $\|A\|_2:=\Ec{\|A\|_\mathrm{op}^2}^{1/2}$ with $\| A\|_\mathrm{op}:= \sup_{\|x\|=1} \|A(x)\|$. Suppose $(X_n)$ obeys (\ref{eq:rec_gen}) and 
\begin{enumerate}
\item[{(A}1{)}] \textsc{Convergence and contraction.}
We have $\|A_r^{(n)}\|_2, \|b_n\|_2<\infty$ for   all $r=1,\ldots,K$ and $n\ge 0$ and there exist random operators $A_1, \ldots, A_K$ on $\Do$ and a $\Do$-valued random variable $b$ with, for some positive sequence $R(n) \downarrow 0$, as $n\to\infty$,
\begin{equation} \label{eq:rate1}
 \|b^{(n)} - b\|_2+\sum_{r=1}^K \left(\|A_r^{(n)} - A_r\|_2 + \left\|{\bf 1}_{\{I^{(n)}_r \le n_0\}} A^{(n)}_r  \right\|_2       \right) = O(R(n))
\end{equation}
and for all $\ell \in \N$,
\begin{equation*}
\E{{\bf 1}_{\{I^{(n)}_r \in \{0,\ldots,\ell\}\cup\{n\}\}} \|A^{(n)}_r\|_\mathrm{op}^2  }\to 0
\end{equation*}
and
\begin{equation} 
L^{*} = \limsup_{n \rightarrow \infty} 
	\E { \sum_{r=1}^K \|A_r^{(n)}\|_\mathrm{op}^2 \frac{R(I_r^{(n)})}{R(n)} } < 1.
\end{equation}
\item[{(A}2{)}] \textsc{Existence and equality of  moments.} $\Ec{\|X_n\|^2} < \infty$ for all $n$ and $\Ec{X_{n_1}(t)} = \Ec{X_{n_2}(t)}$ for all $n_1,n_2 \in \N_0, 
t\in[0,1]$.
\item[{(A}3{)}] \textsc{Existence of a continuous solution.} There exists a solution $X$ of the fixed-point equation \begin{equation} \label{fix}
 X \stackrel{d}{=} \sum_{r=1}^K A_r(X^{(r)}) + b \end{equation} with continuous paths, $\Ec{\|X\|^2} < \infty$ and $\Ec{X(t)} = \Ec{X_1(t)}$ for all $t\in[0,1]$. 
Again $(A_1, \ldots, A_K,b), X^{(1)}, \ldots, X^{(K)}$ are independent and $X^{(1)}, \ldots, X^{(K)}$ are distributed like $X$.
\item[{(A}4{)}] \textsc{Perturbation condition.} $X_n = W_n + h_n$ where $ \|h_n - h \| \to 0$ with $h\in \Do$ and random variables $W_n$ in $\Do$ such that there exists a  sequence $(r_n)$ with, as $n\to\infty$,
\begin{equation*} 
 \Prob{W_n \notin \Drn} \to 0.
\end{equation*}
Here, $\Drn\subset \Do$ denotes the set of functions on the unit interval, for which there is a decomposition of $[0,1]$ into intervals of length as least $r_n$ on which they are constant. 
\item[{(A}5{)}] \textsc{Rate of convergence.} $R(n) =o\left(\log^{-m}(1/r_n)\right)$.
\end{enumerate}

The crucial part that makes everything work consists in choosing a probability metric  in such a way that the limiting map is indeed a contraction. The contraction method presented here for the Banach space $(\Do,\|\,\cdot\,\|)$ is based on the Zolotarev metric $\zeta_s$ and, for our fixed-point equation, we indeed obtain contraction with $s=2$.  This follows by our modified assumption A1 since
$$\E {\sum_{r=1}^K \|A_r\|^2}  = \lim_n \E {\sum_{r=1}^K \|A_r^{(n)} \|^2}\leq \limsup_n \E { \sum_{r=1}^K \|A_r^{(n)}\|^2 \frac{R(I_r^{(n)})}{R(n)} } < 1.$$
The amounts of details to be verified prevents us to provide a complete proof of all the assumptions in the present case. In the remainder of the section, we will not come back on the method and Proposition~\ref{prop:cont} itself but show how it can be applied; we will however, discuss and outline the proof of the main assumptions (A1), (A2), (A3) and (A5).

\begin{prop} \label{prop:cont}
 Let $X_n$ fulfill (\ref{eq:rec_gen}). Provided that Assumptions (A1)--(A3) are satisfied, the solution $X$ of the fixed-point equation (\ref{fix}) is unique.
\begin{compactenum}[i.]
\item For all $t\in[0,1]$, $X_n(t) \to X(t)$ in distribution, with convergence of the first two moments;
\item If $U$ is independent of $(X_n), X$ and distributed on 
$[0,1]$ then $X_n(U)\to X(U)$ in distribution again with convergence of the first two moments. 
\item If  also (A4) and (A5) hold, then $X_n \rightarrow X$ in distribution in $(\Do, \|\cdot\|)$. 
\end{compactenum}
\end{prop}

\subsection{Existence of a continuous solution}\label{sec:limit}

In this section, we outline the proof of existence of a continuous process $Z$ that satisfies the distributional fixed point equation \eqref{eq:limit_process} as it is needed for assumption (A3). 
We construct the process $Z$ as the pointwise limit of martingales. We then show that the convergence is actually almost surely uniform, which allows us 
to conclude that $Z$ is actually continuous with probability one. Write $\mathcal C[0,1]$ for the space of continuous functions on $[0,1]$.

Consider the infinite $4$-ary tree $\mathcal T = \bigcup_{n\ge 0}\{1,2,3,4\}^n$. For a node $u\in \mathcal T$, we write $|u|$ for its depth, i.e.\ the distance between $u$ and the root 
$\varnothing$. The descendants of $u\in \mathcal T$ correspond to all the words in $\mathcal T$ with prefix $u$. Let $\{U_v, v\in \mathcal T\}$ 
and $\{V_v, v\in \mathcal T\}$ be two independent families of i.i.d.\ $[0,1]$-uniform random variables. 

\medskip
\noi\textsc{Construction by iteration.} Define the operator $G:(0,1)^2 \times \mathcal C[0,1]^4 \to \mathcal C[0,1]$ by  
\begin{align}\label{eq:def_G}
G(x,y,f_1,f_2,f_3,f_4)(s)=
&\I{s<x} \left[(xy)^\beta f_1\pran{\frac sx}+ (x(1-y))^\beta f_2\pran{\frac sx}\right]\\
&+\I{s\ge x} \left[((1-x)y)^\beta f_3\pran{\frac {s-x}{1-x}}+ ((1-x)(1-y))^\beta f_4\pran{\frac {s-x}{1-x}}\right].\nonumber
\end{align}
Let $h$ be the map defined by $h(s)=(s(1-s))^{\beta/2}$, where $2\beta=\sqrt{17}-3$. For every node $u\in \mathcal T$, let $Z_0^{u}=h$. Then define 
recursively 
\begin{align} \label{def:recZ}
Z_{n+1}^u=G(U_u, V_u, Z_n^{u1}, Z_n^{u2}, Z_n^{u3}, Z_n^{u4}).
\end{align}
Finally, define $Z_n=Z_n^\varnothing$ to be the value observed at the root of $\mathcal T$ when the iteration has been started with $h$ in all the nodes 
at level $n$. 
%
%

\medskip
\noi\textsc{A series representation for $Z_n$.} For $s\in [0,1]$, $Z_n(s)$ is the sum of exactly $2^n$ terms, each one being the 
contribution of one of the boxes at level $n$ that is cut by the line at $s$. Let $\{Q_i^n(s), 1\le i\le 2^n\}$ be the set of rectangles at level $n$ 
whose first coordinate intersect $s$. Suppose that the projection of $Q_i^n(s)$ on the first coordinate yields the interval $[\ell_i^n, r_i^n]$. Then
\begin{equation} \label{Z_n:explicit}
Z_n(s)=\sum_{i=1}^{2^n} \vol(Q_i^n(s))^\beta \cdot h\pran{\frac{s-\ell_i^n}{r_i^n-\ell_i^n}},
\end{equation}
where $\vol(Q_i^n(s))$ denotes the volume of the rectangle $Q_i^n(s)$. The difference between $Z_n$ and $Z_{n+1}$ only relies in the functions appearing 
the  boxes 
$Q_i^n(s)$: We have \begin{align}\label{eq:Z_n-telescoping}
Z_{n+1}(s)-Z_n(s)=\sum_{i=1}^{2^n} \vol(Q_i^n(s))^\beta \cdot \left[G(U_i',V_i',h,h,h,h)\pran{\frac{s-\ell_i^n}{r_i^n-\ell_i^n}}-h\pran{\frac{s-\ell_i^n}{r_i^n-\ell_i^n}}\right],\end{align}
where $U'_i,V'_i$, $1\le i\le 2^n$ are i.i.d.\ $[0,1]$-uniform random variables. In fact, $U_i'$ and $V_i'$ are some of the variables $U_u, V_u$ for nodes $u$ at
 level $n$. Observe that, although $Q_i^n(s)$ is \emph{not} a product of $n$ independent terms of the form $U V$ because of size-biasing, $U'_i, V'_i$ are in 
fact \emph{unbiased}, i.e. uniform. Let $\mathscr F_n$ denote the $\sigma$-algebra generated by $\{U_u, V_u: |u|< n\}$. 
Then the family $\{U_i',V_i': 1\le i\le 2^n\}$ is independent of $\mathscr F_n$.

\medskip
\noi\textsc{A martingale.} Let $s\in[0,1]$ be fixed. We show that the sequence $(Z_n(s), n\ge 0)$ is a non-negative discrete time martingale ; so it converges with 
probability one to a finite limit $Z(s)$. To prove that $Z_n(s)$ is a indeed a martingale, it suffices to prove that, for $1\le i\le 2^n$,
\[\CE{G(U_i',V_i',h,h,h,h)\pran{\frac{s-\ell_i^n}{r_i^n-\ell_i^n}}}{\mathscr F_n}=h\pran{\frac{s-\ell_i^n}{r_i^n-\ell_i^n}}.\]
Since $U_i',V_i', 1\le i\le 2^n$ are independent of $\mathscr F_n$, this clearly reduces to the following lemma.
\begin{lem}\label{lem:martingale}For the operator $G$ defined in \eqref{eq:def_G} and $U,V$ two independent $[0,1]$-uniform random variables, and any $s\in [0,1]$, we have $\E{G(U,V,h,h,h,h)(s)}=h(s)$.
\end{lem}

\medskip
\noi\textsc{Almost sure continuity.}
%
Assume for the moment that there exist constants $a,b\in (0,1)$ and $C$ such that 
\begin{equation}\label{eq:X_borel-cantelli}
\p{\sup_{s\in [0,1]}|Z_{n+1}(s)-Z_n(s)|\ge a^n} \le C\cdot b^n.
\end{equation}
Then, by the Borel--Cantelli lemma, the sequence $(Z_n)$ is almost surely cauchy with respect to the supremum norm. Completeness of 
$\left( \mathcal C[0,1], \|\cdot \|\right)$ yields the existence of a random process $Z$ with continuous paths such that $Z_n \rightarrow  Z$ uniformly on
$[0,1]$. 
We now move on to showing that there exist constants $a$ and $b$ such that (\ref{eq:X_borel-cantelli}) is satisfied.
We start by a bound for a fixed value $s\in [0,1]$. 

\begin{lem}\label{lem:bound_fixed_s}For every $s\in [0,1]$, any $a\in(0,1)$, and any integer $n$ large enough, we have the bound
$$\p{|Z_{n+1}(s)-Z_n(s)|\ge a^n}\le 4 (16e \log(1/a))^n.$$
\end{lem}

Then, in order to handle the supremum over $s\in [0,1]$, in (\ref{eq:X_borel-cantelli}) note that the number of values taken by $Z_n$ is at most the number of 
boxes at level $n$, i.e.\ $4^n$. 
To avoid unnecessary technicalities, we use fixed points (much more than $4^n$) to control the extent of $\sup_{s \in [0,1]} |Z_{n+1}(s)-Z_n(s)|$. 
Consider the set $V_n$ of $x$-coordinates of the vertical boundaries of all the rectangles at level $n$. Let $L_n=\inf\{|x-y|: x,y\in V_n\}$. Then, on the event 
that $L_n\ge\gamma^n$, we have $$\sup_{s\in [0,1]} |Z_{n+1}(s)-Z_n(s)| \le \sup_{1\le i\le \fl{\gamma^{-n}}} |Z_{n+1}(i \gamma^n)-Z_n(i\gamma^n)|.$$
In particular, it follows by the union bound that, for any $\gamma\in(0,1)$,
\begin{equation*}\label{eq:sup_diff_X}
\p{\sup_{s\in [0,1]} |Z_{n+1}(s)-Z_n(s)|\ge a^n} \le \gamma^{-n}\sup_{s\in [0,1]} \p{|Z_{n+1}(s)-Z_n(s)|\ge a^n} +\p{L_n<\gamma^n}.
\end{equation*}
The following lemma then yields \eqref{eq:X_borel-cantelli} which completes the proof.
\begin{lem}\label{lem:closest_pair}For any positive real number $\gamma$ small enough, it exists an integer $n_1(\gamma)$ with
$$\p{L_n< \gamma^n} \le  6 \cdot 4^{n} \gamma^{n/201}, \qquad n \geq n_1(\gamma).$$
\end{lem}

%
%
\subsection{Uniform convergence of the mean}

The proof of Theorem~\ref{thm:process} requires to show uniform convergence of the first moment $n^{-\beta}\E{C_n(s)}$ towards $\mu_1(s)=K_1(s(1-s))^{\beta/2}$ uniformly on $[0,1]$ in order to verify assumption (A1), in particular the rate $R(n)$ in (\ref{eq:rate1}). Note that, since $C_n(s)$ is continuous in any fixed $s$ almost surely, the function $s \to \E{C_n(s)}$ is continuous for any $n$. 
\citet{CuJo2010} only show pointwise convergence, and proving uniform convergence requires a good deal of additional arguments. 

The first step is to prove a Poissonized version, the fixed-$n$ version is then obtained by a routine Tauberian argument. Consider a Poisson point process with unit intensity on $[0,1]^2\times [0,\infty)$. The first two coordinates represent the location inside the unit square; the third one represents the time of arrival of the point. Let $P_t(s)$ denote the partial match cost for a query at $x=s$ in the quad tree built from the points arrived by time $t$.

\begin{prop}\label{prop:uniform_poisson}There exists $\varepsilon>0$ such that
 $$ \sup_{s \in [0,1]} |t^{-\beta}\Ec{P_t(s)} - \mu_1(s)| = O(t^{-\varepsilon}).$$
\end{prop}

The proof of Proposition~\ref{prop:uniform_poisson} relies crucially on two main ingredients: first, a strengthening of the arguments developed by 
\citet{CuJo2010}, and the speed of convergence $\Ec{C_n(\xi)}$ to $\Ec{\mu_1(\xi)}$ for a uniform query line $\xi$, see \eqref{hwangquad}, by \citet{ChHw2003}. By symmetry, we write for any $\delta\in(0,1/2)$
\begin{align}\label{eq:uniform_decomp}
\sup_{s \in [0,1]} |t^{-\beta}\Ec{P_t(s)} - \mu_1(s)| 
&\le \sup_{s\le \delta} \big|t^{-\beta}\Ec{P_t(s)}-\mu_1(s)\big|+\sup_{s\in(\delta,1/2]} \big|t^{-\beta} \Ec{P_t(s)}-\mu_1(s)\big|.
\end{align}
The two terms in the right-hand side above are controlled by the following two lemmas.  

\begin{lem}[Behavior on the edge]\label{lem:uniform_edge}There exists a constant $C_1$ such that
\begin{equation}\label{eq:uniform_edge}
\limsup_{t\to\infty}\sup_{s\le \delta} \big|t^{-\beta}\Ec{P_t(s)}-\mu_1(s)\big| \le C_1 \delta^{\beta/2}.
\end{equation}
\end{lem}

\begin{lem}[Behavior away from the edge]\label{lem:uniform_middle}There exist constants $C_2,C_3,\eta$ with $0  <  \eta < \beta$  and $\gamma\in (0,1)$ such that, for any integer $k$, and real number $\delta\in (0,1/2)$ we have, for any $t>0$,
\begin{align*}
\sup_{s \geq \delta} |t^{-\beta}\Ec{P_t(s)}-\mu_1(s)| \le C_2 \delta^{-1} (1-\gamma)^k + C_3 k 2^k (\beta-\eta)^{-2k} t^{-\eta}.
\end{align*}
\end{lem}

\medskip
\noi\textsc{Behaviour along the edge.} The behaviour away from the edge is rather involved and we do not describe how the bound in Lemma~\ref{lem:uniform_middle} is obtained. To deal with the term for involving the values of $s\in [0,\delta]$, we relate the value $\Ec{P_t(s)}$ to $\Ec{P_t(\delta)}$. Note that the limit first moment $\mu_1(s)=\lim_{n\to\infty}\Ec{P_t(s)}$ is monotonic for $s\in [0,1/2]$. It seems, at least intuitively, that for any fixed real number $t>0$, $\Ec{P_t(s)}$ should also be monotonic for $s\in[0,1/2]$, but we were unable to prove it. The following weaker version is sufficient for our purpose.

\begin{prop}[Almost monotonicity]\label{prop:monotonicity}
For any $s<1/2$ and $\varepsilon\in [0,1-2s)$, we have
$$\Ec{P_t(s)}\le \E{P_{t(1+\varepsilon)}\pran{\frac{s+\varepsilon}{1+\varepsilon}}}.$$
\end{prop}

\section{Second moment and supremum}\label{sec:moments_sup}

In this section, we obtain explicit expressions about the limit, proving that our general approach also turns out to yield effective and computable results. 

\medskip
\noi\textsc{Variance of the cost.} We first focus on the result in Theorem~\ref{thm:moments}. Our main result implies the convergence $n^{-2\beta}\Ec{C_n(s)^2}\to \Ec{Z(s)^2}$. Write $h(s)=\Ec{Z(s)}=(s(1-s))^{\beta/2}$.
Taking second moments in \eqref{eq:limit_process} and writing it as an integral in terms of $\mu_2(s)=\Ec{Z(s)^2}$ yields that we have
the following integral equation, for every $s\in [0,1]$,
\begin{align*}
\mu_2(s)&=\frac 2{2\beta+1}\left\{\int_s^1 x^{2\beta} \mu_2\pran{\frac s x}dx + \int_0^s (1-x)^{2\beta} \mu_2\pran{\frac {1-s}{1-x}}dx\right\}+ 2 \Be(\beta +1, \beta +1)\cdot \frac{h(s)^2}{\beta+1}.
\end{align*}
One easily verifies that 
the function $f$ given by $f(s)=c_2 h(s)^2$ solves the above equation provided that the constant $c_2$ satisfies
$$c_2=\frac{2}{(2\beta+1)(\beta+1)}c_2+ 2 \frac{ \Be(\beta +1, \beta +1)}{\beta+1}\qquad \text{that is}\qquad c_2=  2  \Be(\beta +1, \beta +1) \frac{2\beta+1}{3(1-\beta)},$$
since $\beta^2 = 2- 3 \beta$.
So if we were sure that $\mu_2(s)$ is indeed $c_2 h(s)^2$, we would have by integration $$\Vc{Z(\xi)} = c_2 \Be( \beta +1, \beta +1) - \Be(\beta/2 +1, \beta/2 +1)^2.$$

To complete the proof, it suffices to show that the integral equation satisfied by $\mu_2$ actually admits a unique solution. To this aim, we show that the map $K$ defined below is a contraction for the supremum norm (the details are omitted)
$$Kf(s)=\frac{2}{2\beta+1}\left\{\int_s^1 x^{2\beta} f\pran{\frac sx}dx + \int_0^s (1-x)^{2\beta} f\pran{\frac{1-s}{1-x}}dx\right\}+2  \Be(\beta +1, \beta +1) \frac{[s(1-s)]^\beta}{\beta+1}.$$

\medskip
\noi\textsc{Cost of the worst query.}  The uniform convergence of $n^{-\beta}C_n(\cdot)$ to the process $Z(\cdot)$ directly implies (continuous mapping theorem) the first claim of Theorem \ref{thm:supremum},
\begin{equation}\label{eq:conv_sup}
\frac{S_n}{K_1 n^{\beta}} \convdist S:=\sup_{s\in [0,1]} Z(s).
\end{equation}
The convergence in the Zolotarev metric $\zeta_2$ on which the contraction method is based here, is strong enough to imply convergence of the first two moments of $S_n$ to the corresponding moments of $S$.

\section{Concluding remarks}

The method we exposed here to obtain refined results about the costs of partial match queries in quadtrees also applies to other geometric data structures based on the divide-and-conquer approach. In particular, similar results can be obtained for the $k$-d trees of \citet{Bentley1975} or the relaxed $k$-d trees of \citet{DuEsMa1998}. 

We conclude by mentioning some open questions. The supremum of the process is of great interest since it upperbounds the cost of any query. Can one identify the moments of the supremum $\sup_{s\in[0,1]}Z(s)$ (first and second)? In the course of our proof, we had to construct a continuous solution of the fixed point equation. We prove convergence in distribution, but conjecture that the convergence actually holds almost surely.

{\small
\setlength{\bibsep}{.2em}
\bibliographystyle{abbrvnat}
\bibliography{bib_quadtrees}

\begin{thebibliography}{35}
\providecommand{\natexlab}[1]{#1}
\providecommand{\url}[1]{\texttt{#1}}
\expandafter\ifx\csname urlstyle\endcsname\relax
  \providecommand{\doi}[1]{doi: #1}\else
  \providecommand{\doi}{doi: \begingroup \urlstyle{rm}\Url}\fi

\bibitem[Bentley(1975)]{Bentley1975}
J.~L. Bentley.
\newblock Multidimensional binary search trees used for associative searching.
\newblock \emph{Communication of the ACM}, 18:\penalty0 509--517, 1975.

\bibitem[Billingsley(1999)]{Billingsley1999}
P.~Billingsley.
\newblock \emph{Convergence of Probability Measures}.
\newblock Wiley Series in Probability and Mathematical Statistics. Wiley,
  second edition, 1999.

\bibitem[Chern and Hwang(2003)]{ChHw2003}
H.~Chern and H.~Hwang.
\newblock Partial match queries in random quadtrees.
\newblock \emph{SIAM Journal on Computing}, 32:\penalty0 904--915, 2003.

\bibitem[Curien and Joseph(2011)]{CuJo2010}
N.~Curien and A.~Joseph.
\newblock {Partial match queries in two-dimensional quadtrees: A probabilistic
  approach}.
\newblock \emph{Advances in Applied Probability}, 43:\penalty0 178--194, 2011.

\bibitem[Drmota et~al.(2008)Drmota, Janson, and Neininger]{DrJaNe}
M.~Drmota, S.~Janson, and R.~Neininger.
\newblock A functional limit theorem for the profile of search trees.
\newblock \emph{Ann. Appl. Probab.}, 18\penalty0 (1):\penalty0 288--333, 2008.
\newblock ISSN 1050-5164.

\bibitem[Duch and Mart\'inez(2002)]{DuMa2002a}
A.~Duch and C.~Mart\'inez.
\newblock On the average performance of orthogonal range search in
  multidimensional data structures.
\newblock \emph{Journal of Algorithms}, 44\penalty0 (1):\penalty0 226--245,
  2002.

\bibitem[Duch et~al.(1998)Duch, Estivill-Castro, and Mart\'inez]{DuEsMa1998}
A.~Duch, V.~Estivill-Castro, and C.~Mart\'inez.
\newblock Randomized $k$-dimensional binary search trees.
\newblock In K.-Y. Chwa and O.~Ibarra, editors, \emph{Proc. of the 9th
  International Symposium on Algorithms and Computation (ISAAC'98)}, volume
  1533 of \emph{Lecture Notes in Computer Science}, pages 199--208. Springer
  Verlag, 1998.

\bibitem[Duch et~al.(2010)Duch, Jim\'enez, and Mart\'inez]{DuJiMa2010}
A.~Duch, R.~Jim\'enez, and C.~Mart\'inez.
\newblock {Rank selection in multidimensional data}.
\newblock In A.~L\'opez-Ortiz, editor, \emph{Proceedings of LATIN}, volume 6034
  of \emph{Lecture Notes in Computer Science}, pages 674--685, Berlin, 2010.
  Springer.

\bibitem[Eickmeyer and R{\"u}schendorf(2007)]{EiRu07}
K.~Eickmeyer and L.~R{\"u}schendorf.
\newblock A limit theorem for recursively defined processes in {$L^p$}.
\newblock \emph{Statist. Decisions}, 25\penalty0 (3):\penalty0 217--235, 2007.
\newblock ISSN 0721-2631.

\bibitem[Finkel and Bentley(1974)]{FiBe1974}
R.~A. Finkel and J.~L. Bentley.
\newblock Quad trees, a data structure for retrieval on composite keys.
\newblock \emph{Acta Informatica}, 4:\penalty0 1--19, 1974.

\bibitem[Flajolet and Lafforgue(1994)]{FlLa1994}
P.~Flajolet and T.~Lafforgue.
\newblock {Search costs in quadtrees and singularity perturbation asymptotics}.
\newblock \emph{Discrete and Computational Geometry}, 12:\penalty0 151--175,
  1994.

\bibitem[Flajolet and Puech(1986)]{FlPu1986}
P.~Flajolet and C.~Puech.
\newblock Partial match retrieval of multidimensional data.
\newblock \emph{Jounal of the ACM}, 33\penalty0 (2):\penalty0 371--407, 1986.

\bibitem[Flajolet and Sedgewick(2009)]{FlSe2009}
P.~Flajolet and R.~Sedgewick.
\newblock \emph{Analytic Combinatorics}.
\newblock Cambridge University Press, Cambridge, UK, 2009.

\bibitem[Flajolet et~al.(1993)Flajolet, Gonnet, Puech, and
  Robson]{FlGoPuRo1993}
P.~Flajolet, G.~H. Gonnet, C.~Puech, and J.~M. Robson.
\newblock Analytic variations on quadtrees.
\newblock \emph{Algorithmica}, 10:\penalty0 473--500, 1993.

\bibitem[Flajolet et~al.(1995)Flajolet, Labelle, Laforest, and
  Salvy]{FlLaLaSa1995}
P.~Flajolet, G.~Labelle, L.~Laforest, and B.~Salvy.
\newblock {Hypergeometrics and the cost structure of quadtrees}.
\newblock \emph{Random Structures and Algorithms}, 7:\penalty0 117--144, 1995.

\bibitem[Gr{\"u}bel(2009)]{GR09}
R.~Gr{\"u}bel.
\newblock On the silhouette of binary search trees.
\newblock \emph{Ann. Appl. Probab.}, 19\penalty0 (5):\penalty0 1781--1802,
  2009.
\newblock ISSN 1050-5164.

\bibitem[Ho-Le(1988)]{hole1988}
K.~Ho-Le.
\newblock {Finite element mesh generation methods: a review and
  classification}.
\newblock \emph{Computer-Aided Design}, 20:\penalty0 27--38, 1988.

\bibitem[Knuth(1998)]{Knuth1998}
D.~E. Knuth.
\newblock \emph{The Art of Computer Programming: Sorting and Searching},
  volume~3.
\newblock Addison-Wesley, 2d edition, 1998.

\bibitem[Mahmoud(1992)]{Mahmoud1992a}
H.~Mahmoud.
\newblock \emph{Evolution of Random Search Trees}.
\newblock Wiley, New York, 1992.

\bibitem[Mart\'inez et~al.(2001)Mart\'inez, Panholzer, and
  Prodinger]{MaPaPr2001}
C.~Mart\'inez, A.~Panholzer, and H.~Prodinger.
\newblock Partial match in relaxed multidimensional search trees.
\newblock \emph{Algorithmica}, 29\penalty0 (1--2):\penalty0 181--204, 2001.

\bibitem[Neininger(2000)]{Neininger2000}
R.~Neininger.
\newblock {Asymptotic distributions for partial match queries in K-d trees}.
\newblock \emph{Random Structures and Algorithms}, 17:\penalty0 403--427, 2000.

\bibitem[Neininger(2001)]{Ne01}
R.~Neininger.
\newblock On a multivariate contraction method for random recursive structures
  with applications to {Q}uicksort.
\newblock \emph{Random Structures Algorithms}, 19\penalty0 (3-4):\penalty0
  498--524, 2001.
\newblock ISSN 1042-9832.
\newblock Analysis of algorithms (Krynica Morska, 2000).

\bibitem[Neininger and R\"uschendorf(2001)]{NeRu2001}
R.~Neininger and L.~R\"uschendorf.
\newblock Limit laws for partial match queries in quadtrees.
\newblock \emph{The Annals of Applied Probability}, 11:\penalty0 452--469,
  2001.

\bibitem[Neininger and R{\"u}schendorf(2004{\natexlab{a}})]{NeRu04}
R.~Neininger and L.~R{\"u}schendorf.
\newblock A general limit theorem for recursive algorithms and combinatorial
  structures.
\newblock \emph{Ann. Appl. Probab.}, 14\penalty0 (1):\penalty0 378--418,
  2004{\natexlab{a}}.
\newblock ISSN 1050-5164.

\bibitem[Neininger and R{\"u}schendorf(2004{\natexlab{b}})]{NeRu04b}
R.~Neininger and L.~R{\"u}schendorf.
\newblock On the contraction method with degenerate limit equation.
\newblock \emph{Ann. Probab.}, 32\penalty0 (3B):\penalty0 2838--2856,
  2004{\natexlab{b}}.
\newblock ISSN 0091-1798.

\bibitem[Neininger and Sulzbach(2011)]{NeSu}
R.~Neininger and H.~Sulzbach.
\newblock On a functional contraction method.
\newblock 2011.
\newblock Manuscript in preparation.

\bibitem[Rachev and R{\"u}schendorf(1995)]{RaRu1995}
S.~Rachev and L.~R{\"u}schendorf.
\newblock Probability metrics and recursive algorithms.
\newblock \emph{Advances in Applied Probability}, 27:\penalty0 770--799, 1995.

\bibitem[Rivest(1976)]{Rivest1976}
R.~Rivest.
\newblock {Partial-match retrieval algorithms}.
\newblock \emph{SIAM Journal on Computing}, 5\penalty0 (19--50), 1976.

\bibitem[R{\"o}sler(1991)]{Roesler1991a}
U.~R{\"o}sler.
\newblock {A limit theorem for ''quicksort''}.
\newblock \emph{RAIRO Informatique th{\'e}orique et Applications}, 25:\penalty0
  85--100, 1991.

\bibitem[R{\"o}sler(1992)]{Roesler1992}
U.~R{\"o}sler.
\newblock A fixed point theorem for distributions.
\newblock \emph{Stochastic Processes and their Applications}, 37:\penalty0
  195--214, 1992.

\bibitem[R{\"o}sler(2001)]{Ro01}
U.~R{\"o}sler.
\newblock On the analysis of stochastic divide and conquer algorithms.
\newblock \emph{Algorithmica}, 29\penalty0 (1-2):\penalty0 238--261, 2001.
\newblock ISSN 0178-4617.
\newblock Average-case analysis of algorithms (Princeton, NJ, 1998).

\bibitem[Samet(1990{\natexlab{a}})]{Samet1990}
H.~Samet.
\newblock \emph{The Design and Analysis of Spatial Data Structures}.
\newblock Addison-Wesley, Reading, MA, 1990{\natexlab{a}}.

\bibitem[Samet(1990{\natexlab{b}})]{Samet1990a}
H.~Samet.
\newblock \emph{Applications of Spatial Data Structures: Computer Graphics,
  Image Processing, and GIS}.
\newblock Addison-Wesley, Reading, MA, 1990{\natexlab{b}}.

\bibitem[Samet(2006)]{Samet2006}
H.~Samet.
\newblock \emph{{Foundations of multidimensional and metric data structures}}.
\newblock Morgan Kaufmann, 2006.

\bibitem[Yerry and Shephard(1983)]{YeSh1983a}
M.~Yerry and M.~Shephard.
\newblock {A modified quadtree approach to finite element mesh generation}.
\newblock \emph{IEEE Computer Graphics and Applications}, 3:\penalty0 39--46,
  1983.

\end{thebibliography}
}


\end{document}